\documentclass[12pt,a4paper]{amsart}
\usepackage{amsmath,amssymb,amscd,amsthm,latexsym} 

\theoremstyle{plain}  % default

\newtheorem{thm}{Theorem}[section]
\newtheorem*{thm*}{Theorem}

\newtheorem{prop}[thm]{Proposition}

\theoremstyle{definition}

\theoremstyle{remark}

\numberwithin{equation}{section}

\renewcommand{\leq}{\leqslant}
\renewcommand{\geq}{\geqslant}
\renewcommand{\setminus}{\smallsetminus}
\newcommand{\R}{\mathbb{R}}
\newcommand{\Z}{\mathbb{Z}}
\newcommand{\C}{\mathbb{C}}

\newcommand{\abs}[1]{\lvert#1\rvert}
\newcommand{\norm}[1]{\lVert#1\rVert}
\newcommand{\floor}[1]{\left[#1\right]}

\DeclareMathOperator{\rk}{rk}

\DeclareMathOperator{\Hom}{Hom}

\DeclareMathOperator*{\Coeff}{Coeff}
\DeclareMathOperator{\Pic}{Pic}

\begin{document}
\title{Topology of $U(2,1)$ representation spaces} 
\author{Peter B. Gothen}
\thanks{
  Partially supported by the Funda\c{c}\~ao para a Ci\^encia e a
  Tecnologia (Portugal) through the Centro de Matem\'atica da
  Universidade do Porto, and by the European Commission through the
  Research Training Networks EAGER and EDGE (Contracts Nos.\ 
  HPRN-CT-2000-00099 and HPRN-CT-2000-00101).
}


\date{September 7, 2001}

\begin{abstract}
  The Betti numbers of moduli spaces of representations of a
  universal central extension of a surface group in the groups
  $U(2,1)$ and $SU(2,1)$ are calculated.  The results are obtained
  using the identification of these moduli spaces with
  moduli spaces of Higgs bundles, and Morse theory, following
  Hitchin's programme \cite{hitchin:1987}.  This requires a careful
  analysis of critical submanifolds which turn out to have a
  description using either symmetric products of the surface or moduli
  spaces of Bradlow pairs.
\end{abstract}

\maketitle

\section{Introduction}

Let $X$ be a closed oriented surface of genus $g \geq 2$.  Consider
the universal central extension
\begin{displaymath}
  0 \to \mathbb{Z} \to \Gamma \to \pi_{1} X \to 0,
\end{displaymath}
generated by the standard generators $a_{1},b_{1}, \ldots,
a_{g},b_{g}$ of $\pi_{1} X$ and a central element $J$ with the
relation
\begin{math}
  J = \prod_{i=1}^{g}[a_{i},b_{i}]
\end{math}.
Our object of study is the moduli space of reductive representations
\begin{displaymath}
  \mathcal{M}_{U(2,1)} = \Hom^{+}(\Gamma, U(2,1)) / U(2,1)
\end{displaymath}
of $\Gamma$ in the non-compact Lie group $U(2,1)$.  As is well known
such representations correspond to flat $U(2,1)$-bundles on the
punctured surface $X \setminus \{p\}$ with fixed holonomy around the
puncture given by the image of the central element $J$ in $U(2,1)$
under the representation.  Such bundles extend (as non-flat bundles)
over the puncture and they are topologically classified by their
reduction to the maximal compact subgroup $U(2) \times U(1) \subseteq
U(2,1)$, that is, by a pair of integers $(d_{1},d_{2})$, where $d_{1}$
is the degree of the rank $2$ complex vector bundle given by
projecting $U(2) \times U(1) \to U(2)$ and $d_{2}$ is the degree of
the complex line bundle given by projecting onto $U(1)$.  These
characteristic numbers are subject to the bound
\begin{equation}
  \label{eq:toledo}
  \abs{d_{1} - 2d_{2}} \leq 3g-3;
\end{equation}
this follows from the work of Domic and Toledo
\cite{domic-toledo:1987} and can also be proved using Higgs bundles
(see Xia \cite{xia:1998}).  Furthermore, Xia proved that
the subspaces $\mathcal{M}_{d_{1},d_{2}}$ of representations with
characteristic numbers $(d_{1},d_{2})$ are exactly the connected
components of $\mathcal{M}_{U(2,1)}$.  

In this paper we calculate the Betti numbers of the spaces
$\mathcal{M}_{d_1,d_2}$ in the case when $(d_1+d_2,3)=1$.  (We need to
impose this condition in order to make sure that the moduli space is
non-singular.)  We use the approach via Higgs bundles and Morse theory
introduced by Hitchin \cite{hitchin:1987}: $\mathcal{M}_{d_1,d_2}$ is
homeomorphic to the moduli space of solutions to a set of equations
from gauge theory known as Hitchin's equations and this
space can be identified with an algebro-geometric moduli space of
Higgs bundles of a certain special form.  The point of view of gauge
theory allows one to do Morse theory in the sense of Bott on the
moduli space and the point of view of algebraic geometry permits a
fairly explicit description of the critical submanifolds in terms of
known spaces: the critical submanifolds turn out to be either closely
related to sym\-me\-tric products of the surface or, more interestingly,
to moduli spaces of Bradlow pairs.

The formula for the Betti numbers of $\mathcal{M}_{U(2,1)}$ is given
in Theorem \ref{thm:poincare-U21} and is fairly complicated.  Of
course one can obtain more explicit results in low genus: see
\eqref{eq:poincare-M_01} and \eqref{eq:poincare-M_10} for the
Poincar\'e polynomials of the two connected components of
$\mathcal{M}_{U(2,1)}$ in the case $g=2$ and $d=1$.

A minor modification of our calculations gives the Betti numbers of
the closely related moduli space $\mathcal{M}_{SU(2,1)}$ of reductive
representations of $\Gamma$ in $SU(2,1)$ (Theorem
\ref{thm:poincare-SU21}).  This space has an interpretation as a
moduli space of fixed determinant Higgs bundles and
$\mathcal{M}_{U(2,1)}$ fibres over the Jacobian of $X$ with fibres
isomorphic to $\mathcal{M}_{SU(2,1)}$.  
By analogy to the case of the moduli space of
stable vector bundles one might expect the Poincar\'e polynomial of
the non-fixed determinant moduli space to be the product of that of the
fixed determinant moduli space by that of the Jacobian (see Atiyah and
Bott \cite{atiyah-bott:1982} and Harder and Narasimhan
\cite{harder-narasimhan:1975}).  It is noteworthy that this is not the
case in our situation and, in particular, it follows that the group
of $3$-torsion points in the Jacobian of $X$ acts
non-trivially on the rational cohomology of $\mathcal{M}_{SU(2,1)}$
(Proposition~\ref{prop:non-trivial-action}).

Another interesting aspect is that the Euler characteristic of the
moduli spaces can be calculated.  The components of
$\mathcal{M}_{U(2,1)}$ all have zero Euler characteristic---this of
course already follows from the fact that they fibre over the Jacobian
which itself has zero Euler characteristic.  More interestingly, the
components of $\mathcal{M}_{SU(2,1)}$ have non-zero Euler
characteristic (see \eqref{eq:euler-fixed}).  Again this is in
contrast to the case of the moduli space of stable bundles of fixed
determinant which has zero Euler characteristic.

This paper is organized as follows: in Section \ref{sec:Higgs-Morse}
we recall the necessary background on Higgs bundles and the Morse
theory strategy; in Section \ref{sec:critical-submanifolds} we analyze
the critical submanifolds and determine their Betti numbers; finally,
in Section \ref{sec:fixed-det}, we treat the fixed determinant moduli
spaces.

\section{Higgs bundles and Morse theory on the moduli space}
\label{sec:Higgs-Morse}

In this section we outline the strategy of our calculations and recall
the necessary background.  For details on this material see the papers
of Corlette \cite{corlette:1988}, Donaldson \cite{donaldson:1987},
Hitchin \cite{hitchin:1987,hitchin:1992}, and Simpson
\cite{simpson:1994a,simpson:1994b}.

Give $X$ the structure of a Riemann surface.
The space $\mathcal{M}_{U(2,1)}$ is homeomorphic to the moduli space
of poly-stable Higgs bundles $(E,\phi)$ of the form
\begin{equation}
  \label{eq:u21-higgs-bundle}
  \begin{aligned}
  E &= E_1 \oplus E_2 \\
  \phi &=
  \left(
  \begin{smallmatrix}
    0 & b \\
    c & 0
  \end{smallmatrix}
  \right),
  \end{aligned}
\end{equation}
where $E_{1}$ is a rank $2$ holomorphic bundle and $E_{2}$ is a
holomorphic line bundle on $X$.  Furthermore, the Higgs field $\phi$ 
consists of two holomorphic maps
\begin{displaymath}
  b \colon E_{2} \to E_{1} \otimes K, \quad
  c \colon E_{1} \to E_{2} \otimes K,
\end{displaymath}
where $K$ is the canonical bundle of $X$.  A Higgs bundle $(E,\phi)$
is called stable if the usual slope stability condition $\mu(F) <
\mu(E)$ is satisfied for any proper non-zero $\phi$-invariant
subbundle $F \subseteq E$ (recall that the slope of a holomorphic
bundle $E$ is $\mu(E) = \deg(E)/\rk(E)$).  In fact it is sufficient to
consider subbundles of the form $F = F_{1} \oplus E_{2}$, where $F_{i}
\subseteq E_{i}$, $i=1,2$ (cf.\ \cite{gothen:2000}).  A Higgs bundle
$(E,\phi)$ is said to be poly-stable if it is the direct sum of stable
Higgs bundles, all of the same slope.

It will be convenient to express the bound \eqref{eq:toledo} in terms
of $d_{2}$ and $d=d_{1} + d_{2}$:
\begin{equation}
  d/3 - (g-1) \leq d_{2} \leq d/3 + (g-1);
  \label{eq:d_2-bound}
\end{equation}
for fixed $d$ there is thus one connected component
$\mathcal{M}_{d_{1},d_{2}}$ for each value of $d_{2}$ in this range.

Note that, clearly, $d_{1} = \deg(E_{1})$ and $d_{2} = \deg(E_{2})$. 
For purposes of topology we can therefore identify
$\mathcal{M}_{d_{1},d_{2}}$ with the moduli space of poly-stable Higgs
bundles of the form \eqref{eq:u21-higgs-bundle} with $\deg(E_{i}) =
d_{i}$, $i=1,2$.  Note also that taking a Higgs bundle of the form
\eqref{eq:u21-higgs-bundle} to its dual $(E_{1}^* \oplus E_{2}^*,
\phi^{t})$ gives an isomorphism of the corresponding components of the
moduli space.  We may therefore assume that $\mu(E_{1}) \leq
\mu(E_{2})$ or, equivalently, that $3d_{2} - d \geq 0$.  This,
together with \eqref{eq:d_2-bound}, gives the range
\begin{equation}
  d/3 \leq d_{2} \leq d/3 + g-1
  \label{eq:d_2-range}
\end{equation}
for $d_{2}$.

If $\deg(E) = d = d_{1} + d_{2}$ is co-prime to $3 = 
\rk(E)$ there are no strictly poly-stable Higgs bundles and in this 
case the moduli space is smooth.  This is essential to doing Morse 
theory on it so we shall make this assumption from now on.

Considering the moduli space from the point of view of gauge 
theory allows one to have metrics on the Riemann surface and the 
bundles $E_{1}$ and $E_{2}$.  It therefore makes sense to consider the 
function
\begin{displaymath}
  f = \norm{\phi}^{2}
\end{displaymath}
on the moduli space.  This function is a perfect Morse-Bott function
and so can be used to calculate the Poincar\'e polynomial of the
moduli space:
\begin{align}
    P_{t}(\mathcal{M}_{U(2,1)}) 
    &= \sum_{i} \dim(H^{i}(\mathcal{M}_{U(2,1)};\mathbb{Q}))t^{i} \notag \\
    &= \sum_{N} t^{\lambda_{N}} P_{t}(N)
  \label{eq:poincare-sum}
\end{align}
where the sum is over the critical submanifolds $N$ of $f$, and 
the index $\lambda_{N}$ is the real dimension of the subbundle of 
the normal bundle of $N$ on which the Hessian of $f$ is negative 
definite.

In order to carry out the calculation it is therefore necessary to be 
able do determine the critical submanifolds of $f$ and their indices: 
a Higgs bundle $(E,\phi)$ is a critical point of $f$ if and only if 
it is a variation of Hodge structure, i.e., it is of the form
\begin{displaymath}
  E = F_{1} \oplus \cdots \oplus F_{m},
\end{displaymath}
where the Higgs field $\phi$ maps $F_{i}$ to $F_{i+1} \otimes K$.
Furthermore, in our case each $F_{i}$ must be a subbundle of $E_{1}$
or $E_{2}$.  The Morse indices can be calculated in terms of the
invariants of the bundles $F_{i}$ (see Section 2.5 of
\cite{gothen:2000}): setting
\begin{math}
  U_{k} = \bigoplus_{k = i-j}\Hom(F_{j},F_{i})
\end{math},
the Morse index at the critical point corresponding to $(E,\phi)$ is 
\begin{equation}
  \label{eq:morse-index}
  \lambda = 2 \sum_{k = 2}^{m-1} \bigl((g-1) \rk (U_k) 
  + (-1)^{k+1} \deg(U_k)\bigr).
\end{equation}
In a similar manner the complex dimension of the critical submanifold 
containing $(E,\phi)$ is
\begin{equation}
  \label{eq:dim-critical}
  1 + (g-1)\bigl(\rk(U_{1}) + \rk(U_{0})\bigr) + \deg(U_{1}) -
  \deg(U_{0})
\end{equation}
and so the complex dimension of the downwards Morse flow of the
critical submanifold through $(E,\phi)$ is given by
\begin{equation}
  \label{eq:dim-downwards-morse}
  1 + \sum_{k = 0}^{m-1} \bigl((g-1) \rk (U_k) 
  + (-1)^{k+1} \deg(U_k)\bigr).
\end{equation}
From this and the determination below of the critical submanifolds one
can easily show that the dimension of the downwards Morse flow is not
the same for all critical submanifolds.  This is in contrast to the
case of moduli spaces of representations in a complex group: it was
shown in \cite{gothen:2000} that in this
case the dimension of the downwards Morse flow is exactly half the
dimension of the moduli space reflecting two fundamental facts about
the moduli space of Higgs bundles: Hausel's theorem \cite{hausel:1998}
that the downwards Morse flow coincides with the nilpotent cone (the
fibre over $0$ of the Hitchin map) and Laumon's theorem
\cite{laumon:1988} that the nilpotent cone is a Lagrangian
submanifold.

\section{Critical submanifolds: Bradlow pairs and symmetric products}
\label{sec:critical-submanifolds}

Next we turn to the detailed analysis of the critical submanifolds. 
This is analogous to the analysis in \cite{gothen:1994}, where the
Betti numbers for the moduli space of rank $3$ Higgs bundles
(corresponding to representations of $\Gamma$ in $SL(3,\C)$) were
calculated.

Note that a Higgs bundle of the form \eqref{eq:u21-higgs-bundle} with 
$\phi = 0$ cannot be stable since at least one of the $\phi$-invariant 
subbundles $E_{1}$ and $E_{2}$ will violate the stability condition.  
It follows that a critical point is represented by a chain $E = 
\bigoplus_{i=1}^{m} F_{i}$ of length $m=2$ or $m=3$.  

It turns out that the length $2$ critical points are essentially what
is known as \emph{holomorphic triples}.  These are generalizations of
Bradlow pairs \cite{bradlow:1991} and were introduced by
Garc{\'\i}a-Prada in \cite{garcia-prada:1994}; they
were later studied systematically by Bradlow and Garc{\'\i}a-Prada
in \cite{bradlow-garcia-prada:1996}.  We briefly recall the relevant
definitions: a holomorphic triple
$(\mathcal{E}_{1},\mathcal{E}_{2},\Phi)$ consists of two holomorphic
vector bundles $\mathcal{E}_{1}$ and $\mathcal{E}_{2}$ and a
holomorphic map $\Phi \colon \mathcal{E}_{2} \to \mathcal{E}_{1}$.  A
holomorphic sub-triple is defined in the obvious way.  For $\alpha \in
\R$ the triple $(\mathcal{E}_{1},\mathcal{E}_{2},\Phi)$ is said to be
$\alpha$-stable if
\begin{displaymath}
  \mu(\mathcal{E}_{1}' \oplus \mathcal{E}_{2}') 
   + \alpha\frac{\rk(\mathcal{E}_{2}')}
   {\rk(\mathcal{E}_{1}') + \rk(\mathcal{E}_{2}')} <
  \mu(\mathcal{E}_{1} \oplus \mathcal{E}_{2}) 
   + \alpha\frac{\rk(\mathcal{E}_{2})}
   {\rk(\mathcal{E}_{1}) + \rk(\mathcal{E}_{2})}
\end{displaymath}
for any proper non-trivial sub-triple
$(\mathcal{E}_{1}',\mathcal{E}_{2}',\Phi')$.

We have the
following proposition concerning critical points represented by length
$2$ chains.

\begin{prop}
  There is one critical submanifold $\mathcal{N}^{2}$ of
  $\mathcal{M}_{d_{1},d_{2}}$ consisting of length $2$ chains.  This
  critical submanifold is isomorphic to the moduli space of
  $\alpha$-stable holomorphic triples
  $(\mathcal{E}_{1},\mathcal{E}_{2},\Phi)$ where $\alpha = 2g-2$,
  $\rk(\mathcal{E}_{1})=2$, $\rk(\mathcal{E}_{2})=1$, 
  $\deg(\mathcal{E}_{1}) = 4g-4 + d_{1}$, and $\deg(\mathcal{E}_{2}) = 
  d_{2}$.  Furthermore, the Morse index of $\mathcal{N}^{2}$ is 
  \begin{displaymath}
    \lambda({\mathcal{N}^{2}}) = 0,
  \end{displaymath}
\end{prop}

\begin{proof}
  This is analogous to Proposition 2.9 of \cite{gothen:1994}, (cf.\ 
  also Theorem 5.2 of \cite{gothen:2000}): under our assumptions a
  length $2$ chain must have $F_1 = E_2$, $F_2 = E_1$ and $c=0$; thus
  setting $\mathcal{E}_1 = E_1 \otimes K$, $\mathcal{E}_2 = E_2$ and
  $\Phi = b$ one obtains a holomorphic triple
  $(\mathcal{E}_1,\mathcal{E}_2,\Phi)$.  One then proves that the
  stability conditions coincide.
  
  The Morse index is obviously zero from \eqref{eq:morse-index}.
\end{proof}

It remains to determine the Poincar\'e polynomial of $\mathcal{N}^2$.
As shown by Garc{\'\i}a-Prada in \cite{garcia-prada:1994}, the fact
that $\mathcal{E}_2$ is a line bundle implies that there is an
isomorphism
\begin{align*}
  \mathcal{N}^2 
    &\to \mathcal{M}^{\text{pairs}} \times \Pic^{d_2} X \\
  (\mathcal{E}_1,\mathcal{E}_2,\Phi) 
    &\mapsto 
    \bigl((\mathcal{E}_2^* \otimes \mathcal{E}_1, \Phi),
    \mathcal{E}_2 \bigr),
\end{align*}
where $\mathcal{M}^{\text{pairs}}$ is the moduli space of
$\alpha$-stable Bradlow pairs $(V,\Phi)$.  
Hence $P_t(\mathcal{N}^2)
= (1+t)^{2g} P_t(\mathcal{M}^{\text{pairs}})$.  The Poincar\'e
polynomial of $\mathcal{M}^{\text{pairs}}$ was, essentially,
determined by Thaddeus in \cite{thaddeus:1994}: he considered the
moduli space of fixed determinant pairs, however (cf.\ Bradlow,
Daskalopoulos and Wentworth
\cite{bradlow-daskalopoulos-wentworth:1996}), the arguments go through
in the case of non-fixed determinant pairs.  The result is that
$\mathcal{N}^2$ has Poincar\'e polynomial
\begin{multline}
  \label{eq:poincare-triples}
  P_t(\mathcal{N}^2) =
  \frac{(1+t)^{4g}}{1-t^2} \\
  \cdot \Coeff_{x^i}
  \left(
  \frac{t^{2\deg(V)+2g-2-4i}}{xt^4-1}-\frac{t^{2i+2}}{x-t^2}
  \right)
  \left(
  \frac{(1+xt)^{2g}}{(1-x)(1-xt^2)}
  \right),
\end{multline}
where 
\begin{displaymath}
  i = \floor{\frac{2d}{3}} - 2d_2 + 2g-2
\end{displaymath}
and $V = \mathcal{E}_2^* \otimes \mathcal{E}_1 = E_2^* \otimes E_1
\otimes K$ so that 
\begin{displaymath}
  \deg(V) = 4g-4 + d - 3d_2.
\end{displaymath}

With regard to the critical points represented by length $3$ chains
note that these are necessarily of the form 
\begin{align*}
  E_1 &= F_1 \oplus F_3 \\
  E_2 &= F_2,
\end{align*}
where the $F_i$ are line bundles and $\phi_i \colon F_i \to F_{i+1}
\otimes K$.  Note also that $c = \phi_1 \in H^0(F_1^{-1}F_2K)$ and $b =
\phi_2 \in H^0(F_2^{-1}F_3K)$ and that stability of $(E,\phi)$ implies
that $b$ and $c$ are non-zero.  Denote the critical submanifold of
length $3$ chains $E=F_1 \oplus F_2 \oplus F_3$ with $\deg(F_i) =
\delta_i$, $i=1,2,3$ by $\mathcal{N}^3(\delta_1,\delta_2,\delta_3)$.
Clearly, $E_1 = F_1 \oplus F_3$ and $E_2 = F_2$, in particular
$\delta_2 = d_2$.
With these preliminaries we have the following description of the
length $3$ critical submanifolds.

\begin{prop}
  There is an isomorphism
  \begin{align*}
    \mathcal{N}^3(\delta_1,\delta_2,\delta_3) &\to S^{m_1}X \times
    S^{m_2}X \times \Pic^{\delta_2}(X) \\
    \bigl(F_1 \oplus F_2 \oplus F_3 ,
    \left(
    \begin{smallmatrix}
      0 & b \\ c & 0
    \end{smallmatrix}
    \right)
    \bigr)
    &\mapsto \bigl((c),(b),F_2\bigr),
  \end{align*}
  where $S^m X$ denotes the $m$th symmetric product of $X$ and
  \begin{align*}
    m_1 &= 2g-2 + \delta_2 - \delta_1, \\
    m_2 &= 2g-2 + \delta_3 - \delta_2.
  \end{align*}
  Furthermore, the Morse index of
  $\mathcal{N}^3(\delta_1,\delta_2,\delta_3)$ is
  \begin{displaymath}
    \lambda\bigl(\mathcal{N}^3(\delta_1,\delta_2,\delta_3)\bigr) = 2g-2 +
    2\delta_1 - 2\delta_3.
  \end{displaymath}
\end{prop}

\begin{proof}
  It is clear that $F_2$ and the divisors $(b)$ and $(c)$ determine
  the bundles $F_1$, $F_2$ and $F_3$ and the sections $b$ and $c$ up
  to scalar multiplication.  It is easy to check that any two Higgs
  bundles obtained in this way are isomorphic and hence the map of the
  statement of the proposition is an isomorphism.

  To calculate the Morse index, one simply applies
  \eqref{eq:morse-index}, noting that $U_2 = \Hom(F_1,F_3)$ and so
  $\deg(U) = \delta_3 - \delta_1$ and $\rk(U_2) = 1$.
\end{proof}

The Poincar\'e polynomial of
$\mathcal{N}^3(\delta_1,\delta_2,\delta_3)$ is calculated from
Macdonald's formula \cite{macdonald:1962} for the Poincar\'e
polynomial of the symmetric product of an algebraic curve to be
\begin{multline}
  P_{t}(\mathcal{N}^3(\delta_1,\delta_2,\delta_3)) = \\ (1+t)^{2g}
  \Coeff_{x^{m_{1}}}\frac{(1+xt)^{2g}}{(1-x)(1-xt^2)}
  \Coeff_{x^{m_{2}}}\frac{(1+xt)^{2g}}{(1-x)(1-xt^2)}
  \label{eq:poincare-symmetric},
\end{multline}
where $m_{1}$ and $m_{2}$ were defined above.

This result is, however, not sufficient: for each value of $d_2 =
\delta_2$ in the range \eqref{eq:d_2-range} we also need to determine the
possible values of the invariants $\delta_1$ and $\delta_3$ (or,
equivalently, the invariants $m_1$ and $m_2$).  To do
this, note first that since $m_1$ and $m_2$ are degrees of line
bundles with non-zero sections, we must have 
\begin{align}
  m_1 &\geq 0 \label{eq:m_1-geq-0}, \\
  m_2 &\geq 0 \label{eq:m_2-geq-0}.
\end{align}
However, $m_1 - m_2 = 2\delta_2 - \delta_1 - \delta_3 = 3d_2 - d$
which is strictly positive by \eqref{eq:d_2-range}.  Hence
\eqref{eq:m_2-geq-0} implies \eqref{eq:m_1-geq-0}.  Secondly, we get
from stability applied to the bundles $F_2 \oplus F_3$ and $F_3$ that
\begin{align}
  \delta_2 + \delta_3 &< \frac{2}{3}d \label{eq:d_2d_3<}, \\
  \delta_3 &< \frac{1}{3}d \label{eq:d_3<}.
\end{align}
In this case \eqref{eq:d_2-range} shows that \eqref{eq:d_2d_3<}
implies \eqref{eq:d_3<}.  The former inequality is equivalent to
\begin{equation}
  \label{eq:m_2<}
  m_2 < 2g-2 + \frac{2}{3}d - 2d_2.
\end{equation}
Note that $m_1$ (and hence $\delta_1$ and $\delta_3$) can be recovered
from $m_2$, $d$ and $d_2$:
\begin{displaymath}
  m_1 = m_2 + 3d_2 - d.
\end{displaymath}
It follows that there is a non-empty critical submanifold
\begin{displaymath}
  \mathcal{N}^3(m_2) = \mathcal{N}^3(\delta_1,\delta_2,\delta_3)\bigr)
\end{displaymath}
for each $m_2$ satisfying \eqref{eq:m_2-geq-0} and \eqref{eq:m_2<}.
It remains to express the Morse index in terms of $m_2$; this is a
simple calculation giving
\begin{equation}
  \label{eq:morse-index-length3}
  \lambda\bigl(\mathcal{N}^3(m_2)\bigr) =
  2(5g-5 + d - 3d_2 - 2m_2).
\end{equation}

We now have all the ingredients required for the calculation of the
Poin\-car{\'e} polynomial of $\mathcal{M}_{U(2,1)}$.

\begin{thm}
  \label{thm:poincare-U21}
  Suppose that $(d,3) = 1$.
  The Poincare polynomial of the component $\mathcal{M}_{d_1,d_2}$ of 
  $\mathcal{M}_{U(2,1)}$ is
  \begin{equation}
    P_{t}(\mathcal{M}_{d_1,d_2}) = P_{t}(\mathcal{N}^{2}) + 
    \sum_{m_{2}=0}^{i}t^{2(5g-5 + d - 3d_2 - 2m_2)}
    P_{t}(\mathcal{N}^{3}(m_{2})),
  \end{equation}
  where $d = d_{1} + d_{2}$, $i = \floor{\frac{2d}{3}} - 2d_2 + 2g-2$ 
  and $P_{t}(\mathcal{N}^{2})$ and $P_{t}(\mathcal{N}^{3}(m_{2}))$ 
  are given by \eqref{eq:poincare-triples} and 
  \eqref{eq:poincare-symmetric} respectively.
\end{thm}

It seems difficult to simplify further this expression and obtain a
closed formula for the Poincar\'e polynomial.  On the other hand, we
can obtain explicit formulas in low genus.  For example, consider $g=2$
and $d=1$, then the values allowed by \eqref{eq:d_2-bound} for $d_{2}$
are $d_{2}=0,1$, in particular $\mathcal{M}_{U(2,1)}$ has two
connected components.  In the case $d_2=1$ (and hence $d_1=0$) we
obtain
\begin{multline}
  \label{eq:poincare-M_01}
  P_{t}(\mathcal{M}_{0,1}) =
  t^{14} + 8\,t^{13} + 30\,t^{12} + 68\,t^{11} + 105\,t^{10} + 
  124\,t^{9} \\ + 128\,t^{8} 
  + 128\,t^{7} + 127\,t^{6} + 120\,t^{5} + 99\,t^{4}
  + 64\,t^{3} + 29\,t^{2} + 8\,t + 1
\end{multline}
(these calculations were performed using the computer algebra system
Maple).  In the case when $d=1$ and $d_2=0$ (and so $d_{1}=1$) the
condition $3d_2-d > 0$ is not satisfied, however, as noted above, the
moduli space is isomorphic to the moduli space for $d=2$ and $d_2=1$,
which does satisfy $3d_2-d > 0$.  In this case one obtains
\begin{multline}
  \label{eq:poincare-M_10}
  P_{t}(\mathcal{M}_{1,0}) = 3\,t^{14} + 28\,t^{13} + 115\,t^{12} 
  + 292\,t^{11} + 528\,t^{10} + 728\,t^{9} \\ + 795\,t^{8}
  + 704\,t^{7} + 511\,t^{6} + 308\,t^{5}
  + 161\,t^{4} + 76\,t^{3} + 30\,t^{2} + 8\,t + 1.
\end{multline}
Note in particular that this shows that the two components are not
homeomorphic.

\section{Fixed determinant bundles and Euler characteristic}
\label{sec:fixed-det}

Consider the determinant map $\mathcal{M}_{d_1,d_2} \to \Pic^{d}(X)$
given by
\begin{displaymath}
  (E,\phi) \mapsto \det(E).
\end{displaymath}
Its fibre over a degree $d$ line bundle $\Lambda$ is naturally
isomorphic to the moduli space of Higgs bundles of the form
\eqref{eq:u21-higgs-bundle} with fixed determinant $\Lambda$ and
$\deg(E_i) = d_i$.  We denote this space by
$\widetilde{\mathcal{M}}_{d_1,d_2}$; it is homeomorphic to the moduli
space of reductive representations of $\Gamma$ in $SU(2,1)$ with the
given invariants $d_1$ and $d_2$.  The moduli space
$\mathcal{M}_{SU(2,1)}$ is the union of the spaces
$\widetilde{\mathcal{M}}_{d_1,d_2}$ over the values of $d_{1}$ and
$d_{2}$ such that \eqref{eq:toledo} is satisfied.  

The calculation of the Poincar\'e polynomial of
$\widetilde{\mathcal{M}}_{d_1,d_2}$ proceeds in the same manner as the
calculation for $\mathcal{M}_{d_1,d_2}$; the main difference is that
the critical submanifolds $\widetilde{\mathcal{N}}^{2}$ and
$\widetilde{\mathcal{N}}^{3}(m_{2})$ now become pull-backs of
$3^{2g}$-fold coverings of the Jacobian in the same way as in
Propositions 2.5 and 3.10 of \cite{gothen:1994}, where the relevant
Poincar\'e polynomials were also calculated.  
% The Poincar\'e
% polynomials of the length $2$ critical submanifolds are simply those
% of the non-fixed determinant case divided by $(1+t)^{2g}$ (the
% Poincar\'e polynomial of the Jacobian), while in the case of the
% length $3$ critical submanifolds there is an extra term due to the
% fact that the critical submanifold is a covering of the product
% of symmetric products of the surface.  
The calculation of the Morse
indices is identical to the non-fixed determinant case.  We omit the
details and only state the result.

\begin{thm}
  \label{thm:poincare-SU21}
  Suppose that $(d,3) = 1$.
  The Poincare polynomial of the component
  $\widetilde{\mathcal{M}}_{d_1,d_2}$ of $\mathcal{M}_{SU(2,1)}$ is
  \begin{equation}
    \label{eq:poincare-fixed}
    P_{t}(\widetilde{\mathcal{M}}_{d_1,d_2}) 
    = P_{t}(\widetilde{\mathcal{N}}^{2}) + 
    \sum_{m_{2}=0}^{i}t^{2(5g-5 + d - 3d_2 - 2m_2)}
    P_{t}(\widetilde{\mathcal{N}}^{3}(m_{2})),
  \end{equation}
  where $d = d_{1} + d_{2}$ and $i = \floor{\frac{2d}{3}} - 2d_2 +
  2g-2$.  The Poincar\'e polynomials $P_{t}(\widetilde{\mathcal{N}}^{2})$
  and $P_{t}(\widetilde{\mathcal{N}}^{3}(m_{2}))$ are given by
  \begin{multline}
    \label{eq:poincare-triples-fixed}
    P_t(\widetilde{\mathcal{N}}^2) = \frac{(1+t)^{2g}}{1-t^2} \\
    \cdot \Coeff_{x^i}
    \left(
    \frac{t^{10g-10 + 2d - 6d_2 -4i}}{xt^4-1}-\frac{t^{2i+2}}{x-t^2}
    \right)
    \left(
    \frac{(1+xt)^{2g}}{(1-x)(1-xt^2)}
    \right),
  \end{multline}
  where $i = \floor{\frac{2d}{3}} - 2d_2 + 2g-2$, and
  \begin{multline}
    P_{t}(\widetilde{\mathcal{N}}^3(m_{2})) = 
    \Coeff_{x^{m_{1}}}\frac{(1+xt)^{2g}}{(1-x)(1-xt^2)}
    \Coeff_{x^{m_{2}}}\frac{(1+xt)^{2g}}{(1-x)(1-xt^2)} \\
    + \binom{2g-2}{m_{1}} \binom{2g-2}{m_{2}}(3^{2g}-1)t^{m_{1}+m_{2}},
    \label{eq:poincare-symmetric-fixed}
  \end{multline}
  where $m_1 = m_2 + 3d_2 - d$.
\end{thm}

It is interesting to note that the Poincar\'e polynomial of
$\mathcal{M}_{d_1,d_2}$ is not simply the product of those of the
Jacobian and $\widetilde{\mathcal{M}}_{d_1,d_2}$, in contrast to the
situation for moduli spaces of stable bundles.
Tensoring with a degree zero line bundle gives an action of the
Jacobian of $X$ on $\mathcal{M}_{U(2,1)}$.  Furthermore, the
determinant map is equivariant if we let the Jacobian act on
$\Pic^{d}(X)$ by tensoring with the third power of a line bundle and
hence
\begin{displaymath}
  \mathcal{M}_{d_1,d_2} \cong
  \bigl(\widetilde{\mathcal{M}}_{d_1,d_2} \times
  \Pic^{d}(X)\bigr) / T_3,
\end{displaymath}
where $T_3 \cong (\Z/3)^{2g}$ is the subgroup of $3$-torsion points of
the Jacobian.  So far this is completely analogous to the case of
moduli of stable vector bundles; see e.g.\ \S 9 of Atiyah and Bott
\cite{atiyah-bott:1982}.  
In that case, the finite covering group acts trivially on the rational
cohomology of the fixed determinant moduli space (this result was
first proved by Harder and Narasimhan \cite{harder-narasimhan:1975})
implying that the Poincar\'e polynomial of the non-fixed determinant
moduli space is simply the product of that of the fixed determinant
moduli space by that of the Jacobian.  Hence our calculations imply
the following result.
\begin{prop}
  \label{prop:non-trivial-action}
  The action of $T_3$ on the rational cohomology of
  $\widetilde{\mathcal{M}}_{d_1,d_2}$ is non-trivial.
\end{prop}
This phenomenon also occurs for moduli of representations of $\Gamma$
in $SL(3,\C)$ (see \cite{gothen:1994}).  
From the point of view of the
Morse theory computation the reason for this result is as follows.
The critical submanifolds in the non-fixed determinant case fibre over
the Jacobian via the determinant map and the fibres are isomorphic to
the fixed determinant critical submanifolds.  One sees from
\cite{gothen:1994} that $T_3$ acts trivially on the the rational
cohomology of the length $2$ critical submanifolds while, on the other
hand, the Poincar\'e polynomial of the length $3$ submanifolds is not
the product of those of the fixed determinant critical submanifold by
that of the Jacobian.  Thus the explanation for
Proposition~\ref{prop:non-trivial-action} from this point of view is
that the action of $T_3$ is non-trivial on the rational cohomology of
the fixed determinant length $3$ submanifolds.

Another point worthy of note is that we can determine the Euler
characteristic of $\widetilde{\mathcal{M}}_{d_1,d_2}$:
from \eqref{eq:poincare-sum} and the fact that the Morse indices are
even it follows that the Euler characteristic of
$\widetilde{\mathcal{M}}_{d_1,d_2}$ is simply the sum of the Euler
characteristics of the critical submanifolds.  The formula
\eqref{eq:poincare-triples-fixed} shows that the critical submanifold
$\widetilde{\mathcal{N}}^2$ has Euler characteristic zero and, since
\begin{displaymath}
  \Coeff_{x^{m_{i}}}\frac{(1+xt)^{2g}}{(1-x)(1-xt^2)}
  = (-1)^{m_{i}}\binom{2g-2}{m_{i}},
\end{displaymath} 
\eqref{eq:poincare-symmetric-fixed} shows that
\begin{align*}
  \chi(\widetilde{\mathcal{N}}^3(m_{2})) &= 
  P_{t}(\widetilde{\mathcal{N}}^3(m_{2}))_{|t=-1} \\
  &= \binom{2g-2}{m_{1}} \binom{2g-2}{m_{2}}
  3^{2g}(-1)^{m_{1}+m_{2}}.
\end{align*}
(This can also be seen directly from the fact that 
$\widetilde{\mathcal{N}}^3(m_{2})$ is a $3^{2g}$-fold covering of 
$S^{m_{1}}X \times S^{m_{2}} X$.)
Noting that $(-1)^{m_{1}+m_{2}} = (-1)^{d+d_{2}}$ we therefore get 
\begin{equation}
    \chi(\widetilde{\mathcal{M}}_{d_1,d_2}) = 
    3^{2g}(-1)^{d+d_{2}}
    \sum_{m_{2}=0}^{i} 
    \binom{2g-2}{m_2 + 3d_2 - d} \binom{2g-2}{m_{2}}.
  \label{eq:euler-fixed}
\end{equation}
Thus, for example, in the case $g=2$, $d=1$, we get
\begin{displaymath}
  \chi(\widetilde{\mathcal{M}}_{0,1}) = 81,
  \quad \chi(\widetilde{\mathcal{M}}_{1,0}) = -324.
\end{displaymath}
In general, we see that $\widetilde{\mathcal{M}}_{d_1,d_2}$ has
non-zero Euler characteristic.  This also happens for representations
in $SL(3,\C)$ (see \cite{gothen:1994}) while it contrasts with the
case of the moduli space of stable bundles with fixed determinant
which has zero Euler characteristic, as one easily sees from the
results of Desale and Ramanan \cite{desale-ramanan:1975}.

\emph{Acknowledgments.} I would like to thank Ivan Smith for asking
the question which prompted this investigation, and Oscar Garc{\'\i}a-Prada
and the referee for useful comments.  I benefitted from participating in
workshops organized by the Vector Bundles on Algebraic Curves group of
the European Research Training Network EAGER.

%\enlargethispage{\baselineskip}
\vspace{1cm}

\noindent
  Departamento de Matem{\'a}tica Pura \\
  Faculdade de Ci{\^e}ncias da Universidade do Porto \\
  4099-002 Porto \\
  Portugal \\
  E-mail: \texttt{pbgothen@fc.up.pt}


\begin{thebibliography}{10}
\newcommand{\enquote}[1]{`#1'}

\bibitem{atiyah-bott:1982}
\textsc{M.~F. Atiyah} and \textsc{R.~Bott}, \enquote{The {Y}ang-{M}ills
  equations over {R}iemann surfaces}, \textit{Philos. Trans. Roy. Soc. London
  Ser. A} 308 (1982) 523--615.

\bibitem{bradlow:1991}
\textsc{S.~B. Bradlow}, \enquote{Special metrics and stability for holomorphic
  bundles with global sections}, \textit{J. Differential Geom.} 33 (1991)
  169--213.

\bibitem{bradlow-daskalopoulos-wentworth:1996}
\textsc{S.~B. Bradlow}, \textsc{G.~D. Daskalopoulos}, and \textsc{R.~A.
  Wentworth}, \enquote{Birational equivalences of vortex moduli},
  \textit{Topology} 35 (1996) 731--748.

\bibitem{bradlow-garcia-prada:1996}
\textsc{S.~B. Bradlow} and \textsc{O.~Garc{\'\i}a-Prada}, \enquote{Stable
  triples, equivariant bundles and dimensional reduction}, \textit{Math. Ann.}
  304 (1996) 225--252.

\bibitem{corlette:1988}
\textsc{K.~Corlette}, \enquote{Flat ${G}$-bundles with canonical metrics},
  \textit{J. Differential Geom.} 28 (1988) 361--382.

\bibitem{desale-ramanan:1975}
\textsc{U.~V. Desale} and \textsc{S.~Ramanan}, \enquote{{P}oincar\'{e}
  polynomials of the variety of stable bundles}, \textit{Math. Ann.} 216 (1975)
  233--244.

\bibitem{domic-toledo:1987}
\textsc{A.~Domic} and \textsc{D.~Toledo}, \enquote{The {G}romov norm of the
  {K}aehler class of symmetric domains}, \textit{Math. Ann.} 276 (1987)
  425--432.

\bibitem{donaldson:1987}
\textsc{S.~K. Donaldson}, \enquote{Twisted harmonic maps and the self-duality
  equations}, \textit{Proc. London Math. Soc. (3)} 55 (1987) 127--131.

\bibitem{garcia-prada:1994}
\textsc{O.~Garc{\'\i}a-Prada}, \enquote{Dimensional reduction of stable
  bundles, vortices and stable pairs}, \textit{Int. J. Math.} 5 (1994) 1--52.

\bibitem{gothen:2000}
\textsc{P.~B. Gothen}, \enquote{Components of spaces of representations and
  stable triples}, Topology, to appear.

\bibitem{gothen:1994}
\textsc{P.~B. Gothen}, \enquote{The {B}etti numbers of the moduli space of
  stable rank $3$ {H}iggs bundles on a {R}iemann surface}, \textit{Int. J.
  Math.} 5 (1994) 861--875.

\bibitem{harder-narasimhan:1975}
\textsc{G.~Harder} and \textsc{M.~S. Narasimhan}, \enquote{On the cohomology
  groups of moduli spaces of vector bundles on curves}, \textit{Math. Ann.} 212
  (1975) 215--248.

\bibitem{hausel:1998}
\textsc{T.~Hausel}, \enquote{Compactification of moduli of {H}iggs bundles},
  \textit{J. Reine Angew. Math.} 503 (1998) 169--192.

\bibitem{hitchin:1987}
\textsc{N.~J. Hitchin}, \enquote{The self-duality equations on a {R}iemann
  surface}, \textit{Proc. London Math. Soc. (3)} 55 (1987) 59--126.

\bibitem{hitchin:1992}
\textsc{N.~J. Hitchin}, \enquote{{L}ie groups and {T}eichm\"{u}ller space},
  \textit{Topology} 31 (1992) 449--473.

\bibitem{laumon:1988}
\textsc{G.~Laumon}, \enquote{Un analogue global du cone nilpotent},
  \textit{Duke Math. J.} 57 (1988) 647--671.

\bibitem{macdonald:1962}
\textsc{I.~G. Macdonald}, \enquote{Symmetric products of an algebraic curve},
  \textit{Topology} 1 (1962) 319--343.

\bibitem{simpson:1994a}
\textsc{C.~T. Simpson}, \enquote{Moduli of representations of the fundamental
  group of a smooth projective variety {I}}, \textit{Inst. Hautes {\'E}tudes
  Sci. Publ. Math.} 79 (1994) 867--918.

\bibitem{simpson:1994b}
\textsc{C.~T. Simpson}, \enquote{Moduli of representations of the fundamental
  group of a smooth projective variety {II}}, \textit{Inst. Hautes {\'E}tudes
  Sci. Publ. Math.} 80 (1994) 5--79.

\bibitem{thaddeus:1994}
\textsc{M.~Thaddeus}, \enquote{Stable pairs, linear systems and the {V}erlinde
  formula}, \textit{Invent. Math.} 117 (1994) 317--353.

\bibitem{xia:1998}
\textsc{E.~Z. Xia}, \enquote{The moduli of flat $\mathrm{PU}(2,1)$ structures
  over {R}iemann surfaces}, \textit{Pacific J. Math.} 195 (2000) 231--256.

\end{thebibliography}
\end{document}